\documentclass[a4paper,11pt,twoside,reqno]{amsart}

\usepackage[utf8]{inputenc}
\usepackage[plainpages=false,pdfpagelabels=true]{hyperref}
\usepackage{amssymb,amsthm}
\usepackage[margin=1in]{geometry}
\usepackage{comment}
\newtheorem{Satz}{Theorem}[section]
\newtheorem{Prop}[Satz]{Proposition}
\newtheorem{Lem}[Satz]{Lemma}
\newtheorem{Thm}[Satz]{Theorem}
\newtheorem{Cor}[Satz]{Corollary}
\theoremstyle{definition}

\newtheorem{Bem}[Satz]{Remark}

\newcommand{\tr}{\operatorname{Tr}}

\newcommand{\C}{{\mathbb{C}}}

\newcommand{\dv}{\text{ }dv}
\newcommand{\s}{\mathbb{S}}

\newcommand{\N}{\mathbb{N}}
\newcommand{\Z}{\mathbb{Z}}
\newcommand{\R}{\mathbb{R}}

\allowdisplaybreaks[1]

\renewcommand{\epsilon}{\varepsilon}

\numberwithin{equation}{section}

\title{On p-harmonic self-maps of spheres}
\author{Volker Branding}
\address{University of Vienna, Faculty of Mathematics\\
Oskar-Morgenstern-Platz 1, 1090 Vienna, Austria\\}
\email{volker.branding@univie.ac.at}

\author{Anna Siffert}
\address{WWU M\"unster, Mathematisches Institut\\
Einsteinstr. 62\\
48149 M\" unster\\
Germany}
\email{anna.siffert@wwu.de}

\date{\today}
\thanks{The first author gratefully acknowledges the support of the Austrian Science Fund (FWF) through the START-project Y963-N35 of Michael Eichmair and
the project ``Geometric analysis of biwave maps'' P 34853.
The second author gratefully acknowledges the
supports of
the Deutsche Forschungsgemeinschaft (DFG, German Research Foundation) - Project-ID 427320536 - SFB 1442, as well as Germany's Excellence Strategy EXC 2044 390685587, Mathematics M\"unster: Dynamics-Geometry-Structure.}


\begin{document}

\begin{abstract}
In this manuscript we study rotationally $p$-harmonic maps between spheres.
We prove that for $p\in\N$ given,
there exist infinitely many $p$-harmonic self-maps of $\s^m$ for each
$m\in\N$ with $p<m< 2+p+2\sqrt{p}$.
In the case of the identity map of \(\s^m\) we explicitly determine the spectrum
of the corresponding Jacobi operator
and show that for \(p>m\), the identity map of $\s^m$ is stable when interpreted
as a \(p\)-harmonic self-map of $\s^m$.
\end{abstract}

\maketitle

\section{Introduction}
Throughout this article let $\phi:(M,g)\rightarrow (N,h)$ be a smooth map between two closed
Riemannian manifolds. We will study \(p\)-harmonic maps, i.e.
critical points of the energy functional
\begin{align}
\label{energy-0}
E_p(\phi)=\frac{1}{p}\int_M|d\phi|^p\dv_g,
\end{align}
where $p\in\N$ with $p\geq 2$.
So far most attention has been paid to the
case $p=2$ in this situation the critical points of \eqref{energy-0}
are precisely \emph{harmonic maps}.

\smallskip

Let us give some motivation for the study of \(p\)-harmonic maps.
The latter are a natural generalization of harmonic maps in the following sense:
It is well-known that the standard harmonic map energy, that is \eqref{energy-0} for \(p=2\),
is invariant under conformal transformations if the domain is two-dimensional.
If the dimension of the domain is larger than two,
this fact no longer holds true.
So one can ask, which energy one could consider in higher dimensions
which is conformally invariant.
It seems that $p$-harmonic maps are the most natural candidate:
It is straightforward to check that the energy \eqref{energy-0} is conformally
invariant in the case that \(p=\dim M\).
Moreover, besides the purely mathematical interest, \(p\)-harmonic maps also have important applications in physics, e.g. in elasticity theory \cite{MR4027845}.

\smallskip

Critical points of $E_p(\phi)$ are characterized as solutions of the
equation
\begin{align}
\label{p-harmonic}
\tr_g\nabla^{\phi^\ast TN}(|d\phi|^{p-2}d\phi)=0.
\end{align}
Note that \eqref{p-harmonic} constitutes a system of quasi-linear elliptic differential equations which degenerates at points where \(|d\phi|=0\).
This additional difficulty leads to the fact that analytic questions on (weak) $p$-harmonic maps require much more
attention compared to the case of standard harmonic maps,
see for example the seminal article of Hardt and Lin \cite{MR896767}.
This comes from the fact that
\eqref{p-harmonic} is not elliptic at points where \(|d\phi|=0\) and one cannot expect
to gain regularity from elliptic theory.
The question of how to deal with weak $p$-harmonic maps was recently taken
up again by Mi\'{s}kiewicz, Petraszczuk and Strzelecki in \cite{mps2022}.
Concerning gradient estimates for \(p\)-harmonic maps we refer to the recent article
of Dong and Lin and references therein \cite{MR4293933}.

\smallskip
Establishing existence results for \(p\)-harmonic maps is a challenging task
and many of the standard techniques from the calculus of variation cannot be applied.
For example, it is not possible to construct a minimizing sequence for \eqref{energy-0}
as is shown in \cite[Section 1]{MR1228587}.

Nevertheless a number of existence results could be achieved over the years:
Both Xu, Yang and Ishida constructed \(p\)-harmonic maps
between spheres by considering the harmonic join of two eigenmaps \cite{MR1416255,MR1228587}.
The existence of \(p\)-harmonic spheres, which are critical points of \eqref{energy-0}
for \(\phi\colon\s^m\to N\), where \(N\) is a compact Riemannian manifold,
was achieved by Kawai, Nakauchi and Takeuchi
in \cite{MR1692995}.
Concerning the heat flow for \(p\)-harmonic maps between compact Riemannian manifolds
existence results were obtained by Fardoun and Regbaoui \cite{MR1948451,MR1951489}.
In \cite{MR1614611} Fardoun proved the existence of rotationally symmetric $p$-harmonic maps, assuming $p\geq 2$, from some Euclidean balls to specific ellipsoids.
Weakly $p$-harmonic maps, for $p\geq 2$, from compact, connected manifolds to closed hemispheres have been studied by Fardoun in \cite{MR2142847}.
The relation between \(p\)-harmonic maps and minimal submanifolds
was studied by Baird and Gudmundsson in
\cite{MR1190447}.
Gastel \cite{MR4027845} established a link between the theory of $p$-harmonic maps and Cosserat models for micropolar elasticity of continua.
An overview on the geometric properties of \(p\)-harmonic maps can be found in the older survey articles \cite{MR1468446,MR1647869}.

\medskip

In this manuscript we study $p$-harmonic self-maps of the $m$-dimensional unit sphere $\s^m$ endowed with the standard round metric. Below we parametrize the sphere $\s^m$ by spherical coordinates $(t,\varphi)$, where $0\leq t\leq\pi$ denotes the colatitude, and $\phi$ gives a point on the equator $\s^{m-1}$ in $\s^m$.
For smooth maps
$$\phi:\s^m\rightarrow\s^m, (t,\varphi)\mapsto (r(t),\varphi)$$
the energy functional (\ref{energy-0}) is given by
\begin{align}
\label{energy-t}
E_p(\phi)=\frac{1}{p}\int_{0}^{\pi}\big(\dot r^2(t)+(m-1)\frac{\sin^2r(t)}{\sin^2t}\big)^\frac{p}{2}\sin^{m-1}tdt.
\end{align}
Each solution of the associated Euler-Lagrange equation
\begin{align}
\label{bvp}
\ddot r(t)&+(m-1)\cot t \dot r(t)-(m-1)\frac{\sin 2r(t)}{2\sin^2t} \\
\nonumber&+(p-2)\dot r(t)\frac{\dot r(t)\ddot r(t)+(m-1)\dot r(t)\frac{\sin 2r(t)}{2\sin^2t}-(m-1)\frac{\sin^ 2 r(t)}{\sin^2t}\cot t}{\dot r^2(t)+(m-1)\frac{\sin^2r(t)}{\sin^2t}}
=0,
\end{align}
with boundary values
\begin{align}
\notag &r(0)=0, \qquad r(\pi)=k\pi, k\in\Z,
\end{align}
for some specific $k$, provides a $p$-harmonic self-map of $\s^m$.

Our main result is the following existence theorem for \(p\)-harmonic self-maps of the Euclidean sphere:

\begin{Thm}
\label{thm:main}
Let $p\in\N$ be given.
For $m\in\N$ with $p<m< 2+p+2\sqrt{p}$ there exist infinitely many $p$-harmonic self-maps of $\s^m$.
\end{Thm}

Our existence result Theorem \ref{thm:main} clearly shows that, despite the analytic difficulties, it can be favorable to
study \(p\)-harmonic maps instead of standard harmonic maps.
In the case \(p=2\), harmonic self-maps of spheres only exist
if \(3\leq m\leq 6\), see the article \cite{MR1436833} by Bizo\'{n} and Chmaj, and Theorem \ref{thm:main} shows that
for \(p>2\) one can find \(p\)-harmonic maps in cases
in which one no longer has an existence result for harmonic maps.

We also show that the assumptions on \(m\) and \(p\) in Theorem \ref{thm:main}
are necessary as is demonstrated by the following non-existence result:

\begin{Thm}
\label{thm:nonexistence}
For \(m>2+p+2\sqrt{p}\) solutions analogous to those constructed in Theorem \ref{thm:main}
do not exist.
\end{Thm}

In the case $p=2$, corresponding to standard harmonic maps,
Theorem \ref{thm:main} and Theorem \ref{thm:nonexistence}
recover the main results of Bizo\'{n} and Chmaj \cite{MR1436833}.
The proof of Theorem \ref{thm:main} makes use of their ideas but due to the additional
nonlinearities we have to overcome a number of substantial technical difficulties.
\medskip

\begin{Bem}
We would like to point out that we could compose the \(\varphi\)-coordinate
with an eigenmap, see for example \cite{MR1228587} for more details,
in which case the factor \((m-1)\) in \eqref{energy-t}
would be shifted such that we would get different conditions for \(m,p\)
in Theorems \ref{thm:main} and \ref{thm:nonexistence}.
\end{Bem}

One important property that characterizes the qualitative behavior of a given $p$-harmonic map is its stability.
If a given $p$-harmonic map is stable, then there does not exist
a second \(p\)-harmonic map \lq nearby\rq, meaning that the critical points of \eqref{energy-0} are isolated.
In this article we study the stability of
the identity map when interpreted
as a \(p\)-harmonic self-map of the sphere.
For more details concerning the stability of
harmonic maps we refer to the introduction of our recent article
\cite{cohom}.

\begin{Thm}
\label{thm:stability-identity}
The eigenvalues of the Jacobi operator describing the stability
of the identity map of \(\s^m\), considered as a \(p\)-harmonic map, are given by
\begin{align}
\lambda_j=-2m+p+j(j+m-1),
\end{align}
where \(j\in\,\mathbb{N}_+\).
Furthermore, the corresponding eigenfunctions can be calculated explicitly,
see Theorem \ref{stability-sphere} below for the details.
\end{Thm}

For \(p=2\), Theorem\,\ref{thm:stability-identity} was already obtained by result of Bizo\'{n} and Chmaj, see \cite[Section 2]{MR1436833}.
From Theorem\,\ref{thm:stability-identity} we obtain the following result.
\begin{Cor}
For \(p>m\), the identity map of $\s^m$ is a stable \(p\)-harmonic self-map of $\s^m$.
\end{Cor}

The stability of the identity map of \(\s^m\) when interpreted as a \(p\)-harmonic map
was already studied in \cite[Section 5]{MR1647869} using abstract methods. In this manuscript we completely solve
the corresponding spectral problem and we also realize that for \(p<m\) the identity map is unstable,
compare \cite[Example 5.2]{MR1647869}.

For a detailed study of the stability of \(p\)-harmonic maps we refer to the articles \cite{MR1197114}
and \cite{MR1145657}.

\medskip

\textbf{Organisation:}
In Section\,\ref{section-existence} we provide the proofs of Theorems\,\ref{thm:main} and \ref{thm:nonexistence}.
After that, Theorem\,\ref{thm:stability-identity} is proved in Section\,\ref{section-stability} and we also give some further details on the second variation of \eqref{energy-0}.

\section{Proofs of Theorems \ref{thm:main} and \ref{thm:nonexistence}}
\label{section-existence}
In this section we provide the proofs of Theorems\,\ref{thm:main} and \ref{thm:nonexistence}.
We may assume $p>2$ since the case $p=2$ has been treated by Bizo\'{n} and Chmaj in \cite{MR1436833}.

\medskip

It turns out to be beneficial to make the following
change of variables
\begin{align*}
t=2\arctan(e^x),\qquad h(x)=r(t)-\frac{\pi}{2}.
\end{align*}

In terms of these variables the energy functional \eqref{energy-t}
acquires the form
\begin{align}
\label{energy-x}
E_p(h)=\frac{1}{p}\int_{-\infty}^{\infty}\big(h'^2(x)+(m-1)\cos^2h(x)\big)^\frac{p}{2}\cosh^{p-m}xdx.
\end{align}

A direct calculation shows that the critical points of \eqref{energy-x} are given by
\begin{align}
\label{euler-lagrange-x-a}
h''(x)+h'(x)(p-m)\tanh x&+\frac{m-1}{2}\sin 2h(x) \\
\nonumber&+h'(x)\frac{p-2}{2}\frac{d}{dx}\log\big(h'^2(x)+(m-1)\cos^2h(x)\big)=0.
\end{align}

Our main tool will be a Lyapunov function $W(x)$
associated with the differential equation \eqref{euler-lagrange-x-a}. To this end, it is convenient to define the function $A(x):\R\rightarrow\R$ by
\begin{align*}
A(x):=h'^2(x)+(m-1)\cos^2h(x)
\end{align*}
such that \eqref{euler-lagrange-x-a} can be written as

\begin{align*}
h''(x)+(p-m)\tanh x\,h'(x)+\frac{m-1}{2}\sin 2h(x)
+\frac{p-2}{2}\frac{A'(x)}{A(x)}h'(x)=0.
\end{align*}

This allows us to derive the following result.

\begin{Lem}
\label{w-monoton}
Let $h(x)$ be a solution of \eqref{euler-lagrange-x-a}.
The function \(W(x):\R\rightarrow\R\) defined by
\begin{align}
\label{lyapunov-a}
W(x):=
A^{\frac{p}{2}-1}(x)\big((p-1)h'^2(x)-(m-1)\cos^2 h(x)\big)
\end{align}
is a Lyapunov function of \eqref{euler-lagrange-x-a}. In particular, $W(x)$ is monotonically increasing
on $[0,\infty)$ if $m\geq p$ and
on $(-\infty,0]$ if $m\leq p$.
\end{Lem}
\begin{proof}
We rewrite \eqref{lyapunov-a} in the following form
\begin{align*}
W(x):=pA^{\tfrac{p}{2}-1}(x)h'^2(x)-A^{\tfrac{p}{2}}(x)
\end{align*}
and calculate
\begin{align*}
W'(x)=&pA^{\tfrac{p}{2}-1}(x)\big(2h'(x)h''(x)
+\frac{p-2}{2}h'^2(x)A^{-1}(x)A'(x)
-\frac{1}{2}A'(x)\big)\\
=&pA^{\frac{p}{2}-1}(x)\big(h'(x)h''(x)+\frac{p-2}{2}h'^2(x)A^{-1}(x)
A'(x)
+\frac{m-1}{2}\sin 2h(x)h'(x)\big) \\
=&p(m-p)A^{\tfrac{p}{2}-1}(x)\tanh x\,h'^2(x),
\end{align*}
whence the claim.
\end{proof}

Note that in the case of \(p=2\) the function \eqref{lyapunov-a} reduces, up to scaling and shifting, to the one considered by Bizo\'{n} and Chmaj in \cite{MR1436833}, see equation (3.8) therein.

\smallskip

For the further analysis it turns out to be useful to rewrite the Euler-Lagrange equation \eqref{euler-lagrange-x-a}
in the following form
\begin{align}
\label{euler-lagrange-x-c}
h''(x)+&(p-m)\tanh x\frac{h'^2(x)+(m-1)\cos^2h(x)}{(p-1)h'^2(x)+(m-1)\cos^2h(x)}h'(x) \\
\nonumber&+\frac{m-1}{2}\frac{(3-p)h'^2(x)+(m-1)\cos^2h(x)}{(p-1)h'^2(x)+(m-1)\cos^2h(x)}\sin 2h(x)=0.
\end{align}

\begin{Bem}
In the case of harmonic self-maps of spheres, that is \(p=2\),
equation \eqref{euler-lagrange-x-c} reduces to
\begin{align}
\label{euler-lagrange-p2}
h''(x)+&(2-m)\tanh xh'(x)+\frac{m-1}{2}\sin 2h(x)=0.
\end{align}
Although at first glance the equation for \(p\)-harmonic self-maps
of spheres looks substantially more complicated than \eqref{euler-lagrange-p2}
we note that both of the fractions in the second and third term
in \eqref{euler-lagrange-x-c} are of magnitude one. Hence, one can expect
that the additional nonlinearities in \eqref{euler-lagrange-x-c}
will lead to additional technical difficulties but can still be handled.
\end{Bem}

From \eqref{euler-lagrange-x-c} we can directly read off that solutions of \eqref{euler-lagrange-x-c} with $h'(0)=0$ are even under the reflection \(x\to-x\) and solutions with $h(0)=0$ are odd under the reflection \(x\to-x\). Below we will focus on such odd and even solutions only. Therefore it is sufficient to consider $x\geq 0$.
An even solution of \eqref{euler-lagrange-x-c} with $h'(0)=b$ will be called \textit{\(b\)-orbit} and denoted by $h(x,b)$.
An odd solution of \eqref{euler-lagrange-x-c} with $h(0)=d$ will be called \textit{\(d\)-orbit}.
In the following we will only pay attention to \(b\)-orbits, but all our considerations can easily be adapted to \(d\)-orbits as well.

We define
\begin{align}
\label{dfn:rho}
\rho^2(x,b):=h^2(x,b)+h'^2(x,b).
\end{align}

The next lemma shows that we have good control over the function \(\rho\)
which can be thought of as the \(W^{1,2}\) norm of \(h\).

\begin{Lem}
Let $h(x)$ be a solution of \eqref{euler-lagrange-x-c}.
Given any \(x_0>0\) and \(\beta>0\), there exists a small number \(\epsilon=\epsilon(\beta,x_0)\)
such that if \(b<\beta\) then \(\rho(x,b)<\beta\) for \(x\leq x_0\).
\end{Lem}

\begin{proof}
We use the shorthand notation $h(x):=h(x,b)$ and $\rho(x):=\rho(x,b)$.
Throughout the proof we make use of \eqref{euler-lagrange-x-c}
as this version of the Euler-Lagrange equation turns out to be most useful
in order to derive estimates.
By assumption \(p>2\) such that
\begin{align*}
\frac{h'^2(x)+(m-1)\cos^2h(x)}{(p-1)h'^2(x)+(m-1)\cos^2h(x)}\leq 1
\end{align*}
and thus we can estimate
\begin{align*}
(m-p)\tanh x\frac{h'^2(x)+(m-1)\cos^2h(x)}{(p-1)h'^2(x)+(m-1)\cos^2h(x)}h'^2(x)\leq
(m-p)h'^2(x).
\end{align*}
Again, since \(p>2\), we have
\begin{align*}
\frac{(3-p)h'^2(x)+(m-1)\cos^2h(x)}{(p-1)h'^2(x)+(m-1)\cos^2h(x)}\leq 1
\end{align*}
such that we may estimate
\begin{align*}
\big|\frac{m-1}{2}&\frac{(3-p)h'^2(x)+(m-1)\cos^2h(x)}{(p-1)h'^2(x)+(m-1)\cos^2h(x)}\sin 2h(x)h'(x)\big| \\
&\leq
\frac{m-1}{2}\big|\sin 2h(x)h'(x)\big|
\leq C |h(x)||h'(x)|.
\end{align*}
Hence, we deduce that
\begin{align*}
\frac{d}{dx}\frac{1}{2}\rho^2(x)=&h(x)h'(x)
+(m-p)\tanh x\frac{h'^2(x)+(m-1)\cos^2h(x)}{(p-1)h'^2(x)+(m-1)\cos^2h(x)}h'^2(x) \\
&-\frac{m-1}{2}\frac{(3-p)h'^2(x)+(m-1)\cos^2h(x)}{(p-1)h'^2(x)+(m-1)\cos^2h(x)}\sin 2h(x)h'(x) \\
\leq &
C\rho^2(x)
\end{align*}
for a positive constant \(C\) which can easily be made explicit.
Integrating the above inequality yields
\begin{align*}
\rho(x)\leq e^{Cx}\rho(0).
\end{align*}
Hence, for any \(x_0>0,\beta>0\) if \(b=\rho(0)<e^{-Cx_0}\beta\) then \(\rho(x)<\beta\) for all \(x\leq x_0\), completing the proof.
\end{proof}

For any \(b\)-orbit we define \(\theta(x,b)\) as follows
\begin{align*}
\theta(0,b):=\frac{\pi}{2},\qquad \theta(x,b):=\arctan\big(\frac{h'(x,b)}{h(x,b)}\big)
\end{align*}
for any \(x>0\). The rotation number \(\Omega(b)\) of the $b$-orbit is given by the expression
\begin{align*}
\Omega(b)=-\frac{1}{\pi}\big(\theta(x_e(b),b)-\theta(0,b)\big).
\end{align*}
Here $x_e(b)$ denotes the smallest $x>0$ at which the $b$-orbit exits the set
\begin{align*}
\Gamma:=\{(h,h',x)\mid h<\tfrac{\pi}{2}, x>0, (h,h')\neq (0,0)\}.
\end{align*}

\begin{Lem}
Let $h(x)$ be a solution of \eqref{euler-lagrange-x-c}.
Assume that \(\rho(x,b)\) is close to zero and that
\begin{align}
\label{assumption-dimension}
p<m<3p-2+2\sqrt{p(p-1)}.
\end{align}
Then \(\theta'(x,b)\) is uniformly bounded from above by a negative constant.
\end{Lem}
\begin{proof}
Again, we make use of the shorthand notation \(\theta(x)=\theta(x,b)\) and \(h(x)=h(x,b)\).
A straightforward calculation yields
\begin{align*}
\theta'(x)=&\frac{h''(x)h(x)-h'^2(x)}{h^2(x)+h'^2(x)}\\
=&-\frac{h'^2(x)}{h^2(x)+h'^2(x)} \\
&+\frac{h(x)}{h^2(x)+h'^2(x)}
(m-p)\tanh x\frac{h'^2(x)+(m-1)\cos^2h(x)}{(p-1)h'^2(x)+(m-1)\cos^2h(x)}h'(x)\\
&-\frac{h(x)}{h^2(x)+h'^2(x)}
\frac{m-1}{2}\frac{(3-p)h'^2(x)+(m-1)\cos^2h(x)}{(p-1)h'^2(x)+(m-1)\cos^2h(x)}\sin 2h(x),
\end{align*}
where we employed \eqref{euler-lagrange-x-c}.
Using a number of trigonometric identities this can manipulated as
\begin{align*}
\theta'(x)=-&\sin^2\theta(x)
+(m-p)\frac{\sin 2\theta(x)}{2}\tanh x\frac{h'^2(x)+(m-1)\cos^2h(x)}{(p-1)h'^2(x)+(m-1)\cos^2h(x)} \\
&-\frac{m-1}{2}\cos^2\theta(x)\frac{(3-p)h'^2(x)+(m-1)\cos^2h(x)}{(p-1)h'^2(x)+(m-1)\cos^2h(x)}\frac{\sin 2h(x)}{h(x)}.
\end{align*}
We may further manipulate the above expression as follows
\begin{align}
\label{eq:theta-der}
\theta'(x)=&-\frac{m}{2}+(\frac{p}{2}-1)\cos 2\theta(x)+\frac{1}{2}(m-p)
\big(|\sin 2\theta(x)|-\cos 2\theta(x)\big)+\delta_p\\
\nonumber=&-\frac{m}{2}+(p-1-\frac{m}{2})\cos 2\theta(x)+\frac{1}{2}(m-p)|\sin 2\theta(x)|
+\delta_p,
\end{align}
where
\begin{align*}
\delta_p:=&
\frac{m-p}{2}\bigg(\tanh(x)\frac{h'^2(x)+(m-1)\cos^2h(x)}{(p-1)h'^2(x)+(m-1)\cos^2h(x)}
\sin 2\theta(x)
-|\sin 2\theta(x)|\bigg)
\\
&+(m-1)\cos^2(\theta(x))\bigg(1-\frac{1}{2}\frac{(3-p)h'^2(x)+(m-1)\cos^2h(x)}{(p-1)h'^2(x)+(m-1)\cos^2h(x)}\frac{\sin 2h(x)}{h(x)}\bigg).
\end{align*}
Now, we set \(a:=(p-1-\frac{m}{2}),b:=\frac{1}{2}(m-p)\) and consider the function
\begin{align*}
g(\theta):=a\cos 2\theta+b|\sin 2\theta|.
\end{align*}
Note that due to the assumption \eqref{assumption-dimension} we have \(b>0\).
The extremal points of \(g(\theta)\) can be characterized by
\begin{align}
\label{extremal-points}
-a+b\frac{\cos 2\theta}{|\sin 2\theta|}=0.
\end{align}
This can for example be achieved by regularizing the absolute value in
the definition of \(g(\theta)\).
In order to solve \eqref{extremal-points} for \(\theta\)
we have to make a case distinction between \(a>0\) and \(a<0\).
If \(a>0\), then a direct calculation shows that \eqref{extremal-points} vanishes for
\(\theta_0=\frac{1}{2}\operatorname{arccot}(\frac{a}{b})\),
in the case that \(a<0\), we find \(\theta_0=-\frac{1}{2}\operatorname{arccot}(\frac{a}{b})\).

However, it is straightforward to check that in both cases
\begin{align*}
g(\theta_0)
=\sqrt{a^2+b^2}.
\end{align*}

Reinserting the definition of \(a,b\) we get
\begin{align*}
\sqrt{a^2+b^2}
&=\frac{m-p}{\sqrt{2}}\sqrt{1+\frac{1}{2}\frac{(2-p)^2}{(m-p)^2}+\frac{2-p}{m-p}}.
\end{align*}

Now, we may estimate
\begin{align*}
\theta'(x)\leq&-\frac{m}{2}+\frac{1}{2}\sqrt{(m-2p+2)^2+(m-p)^2}+\delta_p .\\
\end{align*}

A direct calculation shows that
\begin{align*}
-\frac{m}{2}+\frac{1}{2}\sqrt{(m-2p+2)^2+(m-p)^2}\leq 0
\end{align*}
if \(m<3p-2+2\sqrt{p(p-1)}\).

By the previous Lemma, for any \(x_0>0\), by choosing \(b\) sufficiently small, we can make the second term on the right-hand side of \(\delta_p\) arbitrarily small for \(x\leq x_0\).
Now, for \(m<3p-2+2\sqrt{p(p-1)}\) the right hand side of the above inequality
is bounded from above by a negative constant which completes the proof.
\end{proof}

Combining the previous two lemmas we obtain the following result.

\begin{Prop}
\label{prop-omega}
Assume that \(p<m<3p-2+2\sqrt{p(p-1)}\). Then for any given \(N>0\)
there exists an \(\epsilon>0\) such that if \(0<b<\epsilon\),
then \(\Omega(b)>N\).
\end{Prop}
\begin{proof}
This follows from integrating the equation \(\theta'(x)\leq -c^2\)
and the definition of \(\Omega(b)\).
\end{proof}

\begin{Lem}
\label{lem:connecting-orbit}
If a $b$-orbit stays in \(\Gamma\) for all \(x\geq 0\)
and if in addition \(h(x,b)\) has a finite number of zeros,
then \(\lim_{x\to\infty}h(x,b)\to\pm\frac{\pi}{2}\) and \(\lim_{x\to\infty}h'(x,b)\to 0\).
\end{Lem}

\begin{proof}
For simplicity, we use the shorthand notation $h(x):=h(x,b)$.
Suppose that \(h(x_0)\) has an extremal point for some \(x_0\in\R\),
that is \(h'(x_0)=0\). Then from \eqref{euler-lagrange-x-c} we get
\begin{align*}
h''(x_0)+\frac{m-1}{2}\sin 2h(x_0)=0
\end{align*}
and we may conclude that in the set \(\Gamma\) the solution \(h(x)\)
cannot have a positive minimum or a negative maximum.
We deduce that for \(x\) sufficiently large the solution \(h(x)\)
must be monotonic as it is not allowed to oscillate around zero
at infinity by assumption. Due to the monotonicity of \(h(x)\) we have
the limit \(\lim_{x\to\infty}h'(x)=0\) and we conclude from \eqref{euler-lagrange-x-c}
that \(\lim_{x\to\infty}h(x)=\pm\frac{\pi}{2}\) or \(\lim_{x\to\infty}h(x)=0\).

To complete the proof it remains to show that the second case, i.e. \(\lim_{x\to\infty}h(x)=0\), does not occur.
To this end we assume that \(x_0\) is the last extremum of \(h(x)\) before
\(x\) approaches infinity. Now consider the Lyapunov function \(W(x)\)
at the point \(x_0\): it satisfies \(W(x_0)>-(m-1)^{\tfrac{p}{2}}\). Indeed, if $h'(x_0)=0$, we would have $h(x)\equiv 0$, contradicting that $h(x)$ has a finite number of zeros.
As \(W(x)\) is monotonically
increasing we would get a contradiction from $\lim_{x\to\infty}h(x)=0$ since then we would have $W(\infty)=-(m-1)^{\tfrac{p}{2}}$. Therefore, this case does not occur.
\end{proof}

\begin{proof}[Proof of Theorem \ref{thm:main}]
We define the set
\begin{align*}
S_1:=\{b\mid \textrm{$b$-orbit exits } \Gamma \textrm{ via } h=\frac{\pi}{2}
\textrm{ with } \Omega(b)\leq\tfrac{1}{2}\}.
\end{align*}
In order to ensure that the set \(S_1\) is not the empty set,
we employ the Lyapunov function \(W(x)\) as considered in Lemma \ref{w-monoton}.
Note that for a $b$-orbit we have \(h(0)=0\) and \(h'(0)=b\).
A direct calculation shows that
\begin{align*}
W(0)=\big(b^2+(m-1)\big)^{\frac{p}{2}-1}\big((p-1) b^2-(m-1)\big).
\end{align*}
By choosing \(b>0\) appropriately we get \(W(0)>0\)
and due to the monotonicity of \(W(x)\) established in \eqref{lyapunov-a}
we then have \(W(x)>0\) for all \(x\in\R\).
Therefore we have \(h'(x)>0\) for all \(x>0\).
Consequently, the $b$-orbit constructed above exits the set \(\Gamma\)
through \(h=\frac{\pi}{2}\) with \(\Omega(b)<\tfrac{1}{2}\).
Set \(b_1:=\inf S_1\). By Proposition \ref{prop-omega} we know that \(b_1>0\).
Note that the \(b_1\)-orbit cannot exit the set \(\Gamma\) via \(h=\frac{\pi}{2}\)
as this would also hold true for any \lq nearby orbit\rq\, with \(b<b_1\)
which would contradict the definition of \(b_1\).
Hence, the \(b_1\)-orbit stays in \(\Gamma\) for all \(x>0\)
and by Lemma \ref{lem:connecting-orbit} we get that \(\Omega(b_1)=\tfrac{1}{2}\).

\smallskip

Next, we construct the second orbit \(S_2\)
which we define by
\begin{align*}
S_2:=\{b\mid \textrm{$b$-orbit exits } \Gamma \textrm{ via } h=\frac{\pi}{2}
\textrm{ with } \Omega(b)\leq \tfrac{3}{2}\}.
\end{align*}
The strategy to complete the proof is to show that \(S_2\) is non-empty and then to proceed inductively.

\smallskip

Note that by definition of $b_1$, for \(b<b_1\) we have $\Omega(b)>\tfrac{1}{2}$.
We will show that for $b<b_1$ still sufficiently close to \(b_1\), we have $\Omega(b)\leq \tfrac{3}{2}$ and hence \(b\in S_2\).
For this purpose, let $x_a>0$ be such that $h'(x_a)=0$ and $0<h(x_a)<\tfrac{\pi}{2}$.
By choosing $b$ sufficiently close to $b_1$, we can make $x_a$ as large as we want. Further, let $x_b$ the smallest $x>x_a$ such that $h(x)=0$. Note that $x_a$ and $x_b$ exist by the definition of $b_1$. Moreover, from \eqref{euler-lagrange-x-a} we have
\begin{align}
\label{mon}
h'(x)<0\qquad\mbox{for all}\qquad x\in(x_a,x_b).
\end{align}
Our goal is to prove that for $b$ appropriately chosen we have
$W(x_b)\geq 0$. Hence, due to \eqref{mon} and the monotonicity of $W(x)$, $h'(x)<0$ for all $x\geq x_a$, providing that $S_2$ is non-empty.
To accomplish the above goal we estimate the difference $W(x)-W(x_a)$ in two ways and combine the resulting estimates.

\smallskip

On the one hand, from the proof of Lemma\,\ref{w-monoton}, we have
$$W'(x)=p(m-p)A^{\tfrac{p}{2}-1}(x)\tanh x\,h'^2(x)$$
and thus we get
\begin{align}
\label{w-1}
W(x)-W(x_a)=\notag&p(m-p)\int_{x_a}^xA^{\tfrac{p}{2}-1}(s)\tanh s \,h'^2(s)\,ds\\
\geq \notag& p(m-p)A^{\tfrac{p}{2}-1}(x_a)\tanh x_a \int_{x_a}^x h'^2(s)\,ds\\
=& p(m-p)(m-1)^{\tfrac{p}{2}-1}\cos^{p-2}h(x_a)\tanh x_a \int_{x_a}^x h'^2(s)\,ds
\end{align}
for $x\in (x_a,x_b)$. Here, we made use of the fact that $h(x)$ cannot have a positive minimum and thus $A^{\tfrac{p}{2}-1}(x)$ is monotonically increasing in $(x_a,x_b)$.

\smallskip

On the other hand, the monotonicity of $W(x)$ yields
\begin{align}
\label{w1}
W(x)-W(x_a)\notag=&\big(h'^2(x)+(m-1)\cos^2 h(x)\big)^{\frac{p}{2}-1}\big((p-1)h'^2(x)-(m-1)\cos^2 h(x)\big)\\&+(m-1)^{\tfrac{p}{2}}\cos^p h(x_a)\geq 0.
\end{align}
Below we use the shorthand notations
$$\epsilon:=(m-1)^{\tfrac{p}{2}}\cos^p h(x_a)$$
and
$$f(x):=(p-1)h'^2(x)-(m-1)\cos^2 h(x).$$
Further, let $x_0\in(x_a,x_b)$ such that $h(x_0)=\tfrac{\pi}{4}$.
For $x\in [x_0,x_b]$ we thus have $\cos^2 h(x)\geq\tfrac{1}{2}$.

\smallskip

Note that if there exists a $\hat{x}\geq0$ such that $h'(\hat{x})^2>\tfrac{m-1}{p-1},$
then $W(\hat{x})>0$. Hence $W(x)>0$ for all $x\geq \hat{x}$ and thus $h'(x)\neq0$ for $x\geq \hat{x}$. Consequently, below we may assume without loss of generality
\begin{align}
\label{hstrich}
  h'(x)^2\leq\tfrac{m-1}{p-1}
\end{align}
for all $x\geq 0$.

Let $x\in [x_0,x_b]$ such that $f(x)\geq 0$.
Thus, from \eqref{w1} and \eqref{hstrich} we have
$$(\tfrac{m-1}{p-1}+m-1)^{\tfrac{p}{2}-1}f(x)\geq
\big(h'^2(x)+(m-1)\cos^2 h(x)\big)^{\frac{p}{2}-1}f(x)\geq -\epsilon.$$
Thus, we find
\begin{align}
\label{f1}
h'^2(x)\geq \tfrac{1}{p-1}(\tfrac{m-1}{2}-(\tfrac{m-1}{p-1}+m-1)^{1-\tfrac{p}{2}}\epsilon).
\end{align}

Let $x\in [x_0,x_b]$ such that $f(x)\leq 0$.
Thus, from \eqref{w1} we get
\begin{align*}
((m-1)\cos^2 h(x_a))^{\frac{p}{2}-1}f(x)\geq \big(h'^2(x)+(m-1)\cos^2 h(x)\big)^{\frac{p}{2}-1}f(x)\geq -\epsilon,
\end{align*}
which yields
\begin{align}
\label{f2}
h'^2(x)\geq \tfrac{1}{p-1}(\tfrac{m-1}{2}-(m-1)\cos^2 h(x_a)).
\end{align}
Combining \eqref{f1} and \eqref{f2} we get
\begin{align*}
h'^2(x)\geq \tfrac{1}{p-1}\big(\tfrac{m-1}{2}-\mbox{max}\big[(m-1)\cos^2 h(x_a),(\tfrac{m-1}{p-1}+m-1)^{1-\tfrac{p}{2}}\epsilon\big]\big)=:c_1
\end{align*}
for all $x\in [x_0,x_b]$.
Note that for $p$ and $m$ given, we can chose $x_a$ sufficiently large such that $c_1$ is positive.

\smallskip

From \eqref{hstrich} we obtain a lower positive bound $c_2$ on $x_b-x_0$, which depends on $p$ and $m$.
Thus, from \eqref{w-1} we get
\begin{align*}
W(x_b)&\geq W(x_a)+p(m-p)(m-1)^{\tfrac{p}{2}-1}\cos^{p-2}h(x_a)\tanh x_a \int_{x_0}^{x_b} h'^2(s)\,ds\\
&\geq -(m-1)^{\tfrac{p}{2}}\cos^p h(x_a)+p(m-p)(m-1)^{\tfrac{p}{2}-1}\cos^{p-2}h(x_a)
\tanh x_a c_1c_2\geq 0.
\end{align*}
For the last estimate we assume that $x_a$ is large enough.

\smallskip

Hence, the monotonicity of $W(x)$ implies $h'(x)<0$ for all $x>x_0$.
Consequently, the set $S_2$ is non-empty and we can proceed as before.
We can now iterate the argument.

This completes the proof of Theorem \ref{thm:main}.

\end{proof}

\begin{proof}[Proof of Theorem \ref{thm:nonexistence}]
First, we study the linearization of \eqref{euler-lagrange-x-c}.
It is straightforward to check that the right-hand side of \eqref{euler-lagrange-x-c} vanishes in the case that \(h(x)=k\frac{\pi}{2}\) for \(k\in\Z\).
In order to linearize the above equation we now calculate \(\frac{d}{ds}\big|_{s=0}h(x)=h_L(x)\)
and evaluate it at the critical point \(h_0=0\). The resulting equation is
\begin{align}
\label{linearization}
h_L''(x)=(m-p)\tanh(x)h'_L(x)+(1-m)h_L(x).
\end{align}

Although it is straightforward to explicitly solve
\eqref{linearization} in terms of orthogonal polynomials we can directly
obtain the desired information by investigating \eqref{linearization} for
large values of \(x\) in which case we need to study
\begin{align}
\label{linearization-large-x}
h_L''(x)=(m-p)h'_L(x)+(1-m)h_L(x).
\end{align}
Making the ansatz \(h_L(x)=e^{\alpha x},\alpha\in\C\), we obtain the algebraic equation
\begin{align*}
\alpha^2-(m-p)\alpha+(m-1)=0
\end{align*}
which has the solutions
\begin{align*}
\alpha_\pm=\frac{1}{2}(m-p)\pm\frac{1}{2}\sqrt{m^2-2m(2+p)+p^2+4}.
\end{align*}
It is easy to see that for \(m<2+p+2\sqrt{p}\) we have that \(\alpha\in\C\)
such that \eqref{linearization-large-x} admits oscillatory
solutions. If \(m>2+p+2\sqrt{p}\) the solutions
of \eqref{linearization-large-x} will grow exponentially which is an obstruction
to the existence of solutions of \eqref{euler-lagrange-x-c}.
\end{proof}

\begin{Bem}
We want to point out that there are two places in the proof
of Theorem \ref{thm:main} where we needed to impose a restriction on the dimension.
In Lemma \ref{prop-omega} we needed to make the assumption \(p<m<3p-2+2\sqrt{p(p-1)}\)
in order to get the desired bound on \(\theta'(x)\).
However, in the proof of the above Theorem we have seen that we
have to impose the stronger restriction $p<m< 2+p+2\sqrt{p}$.
Note that in the case of \(p=2\) both inequalities coincide
such that this phenomena seems to be specific for \(p\)-harmonic maps.
\end{Bem}

\section{Second variation formula and Proof of Theorem \ref{thm:stability-identity}}
\label{section-stability}
In this section we prove Theorem \ref{thm:stability-identity}.
The techniques that we employ here are inspired from our recent
investigation of the (equivariant)
stability of harmonic self-maps of cohomogeneity one manifolds
\cite{cohom}.

\smallskip

In order to prove Theorem \ref{thm:stability-identity} we hence
first calculate the second variation of the energy \eqref{energy-x}.
\begin{Lem}
Suppose that \(h(x)\) is a critical point of \eqref{energy-x}.
Then the second variation of \eqref{energy-x} is given by
\begin{align}
\frac{d^2}{ds^2}\big|_{s=0}E_{p}(h_s)
\nonumber=&-\int_{-\infty}^\infty\xi(x)\big(h'^2(x)+(m-1)\cos^2h(x)\big)^{\frac{p}{2}-1}\\
\nonumber &\times\bigg(\xi''(x)+(p-m)\tanh x \xi'(x)+(m-1)\cos 2h(x)\xi(x) \\
\nonumber &+\xi'(x)\frac{p-2}{2}\frac{d}{dx}\log\big(h'^2(x)+(m-1)\cos^2h(x)\big) \\
\nonumber &+h'(x)(p-2)\frac{d}{dx}\big(\frac{h'(x)\xi'(x)-\frac{m-1}{2}\sin 2h(x)\xi(x)}{h'^2(x)+(m-1)\cos^2h(x)}\big)
\bigg)\frac{dx}{\cosh^{m-p}x},
\end{align}
where \(h_s(x)\) is an one-parameter variation of the map \(h(x)\) that satisfies
\begin{align*}
\frac{d}{ds}\big|_{s=0}h_s(x)=\xi(x).
\end{align*}
\end{Lem}

\begin{proof}
This follows by a direct, but lengthy calculation.
\end{proof}

A $p$-harmonic map is \emph{stable} if and only if the Jacobi operator associated to the second variation of the energy \eqref{energy-0} has only positive eigenvalues. Hence, in order to investigate the stability of \(p\)-harmonic self-maps of the sphere we have to solve the following eigenvalue problem
\begin{align}
\label{ode-second-variation}
\xi''(x)&+(p-m)\tanh x \xi'(x)+(m-1)\cos 2h(x)\xi(x) \\
\nonumber &+\xi'(x)\frac{p-2}{2}\frac{d}{dx}\log\big(h'^2(x)+(m-1)\cos^2r(x)\big) \\
\nonumber &+h'(x)(p-2)\frac{d}{dx}\big(\frac{h'(x)\xi'(x)-\frac{m-1}{2}\sin 2h(x)\xi(x)}{h'^2(x)+(m-1)\cos^2h(x)}\big) \\
\nonumber&+\frac{\lambda}{\big(h'^2(x)+(m-1)\cos^2h(x)\big)^{\frac{p}{2}-1}}\frac{1}{\cosh^px}=0.
\end{align}

Below we focus on the stability of the identity map.

\begin{Lem}
The stability of the identity map of \(\s^m\) considered as a solution of the equation for \(p\)-harmonic
self-maps of the sphere is described by the following spectral problem
\begin{align}
\label{jacobi-identity-spectral}
\xi''(x)&+\xi'(x)(2 -m)\tanh x -\xi(x)(m-1)\big(\tanh^2x-\frac{1}{\cosh^2x}\big) \\
&\nonumber+\frac{p-2}{\cosh x }\frac{d}{dx}\big(\frac{\cosh x}{m}(\xi'(x)-(m-1)\tanh x \xi(x))\big) \\
&\nonumber+\frac{\hat{\lambda}}{\cosh^{2}x}\xi(x)=0,
\end{align}
where
$$\hat{\lambda}=\frac{\lambda}{m^{\frac{p}{2}-1}}.$$
\end{Lem}
\begin{proof}
In the case of the identity map \(h_1(x)=-\frac{\pi}{2}+2\arctan(e^x)\) we have the following identities
\begin{align*}
h_1'^2(x)+(m-1)\cos^2h_1(x)&=\frac{m}{\cosh^2 x},\\
\sin 2h_1(x)&=2\frac{\tanh x}{\cosh x},\\
\cos 2h_1(x)&=\frac{1}{\cosh^2x}-\tanh^2x.
\end{align*}
Moreover, a direct calculation yields
\begin{align*}
\xi'(x)\frac{p-2}{2}\frac{d}{dx}\log&\big(h_1'^2(x)+(m-1)\cos^2h_1(x)\big)
=(2-p)\tanh x \xi'(x), \\
h'(x)(p-2)\frac{d}{dx}\big(&\frac{h_1'(x)\xi'(x)-\frac{m-1}{2}\sin 2h_1(x)\xi(x)}{h_1'^2(x)+(m-1)\cos^2h_1(x)}\big) \\
=&\frac{p-2}{\cosh x }\frac{d}{dx}\big(\frac{\cosh x}{m}(\xi'(x)-(m-1)\tanh x \xi(x))\big).
\end{align*}
Inserting these identities into \eqref{ode-second-variation} completes the proof.
\end{proof}

In order to complete the proof of Theorem \ref{thm:stability-identity}
we recall some facts on the so-called Gegenbauer polynomials.
For more details on this subject we refer to \cite[Chapter 22]{MR0167642}
and the website \cite[Chapter 18]{NIST:DLMF}.

Consider a second order linear ordinary differential equation of the form
\begin{align}
\label{ode-gegenbauer}
(1-x^2)f''(x)-(2\alpha+1)xf'(x)+\lambda_jf(x)=0.
\end{align}
with \(\lambda_j=j(j+2\alpha)\).
Then \eqref{ode-gegenbauer} is solved by the \emph{Gegenbauer polynomials}
\(C_j^{(\alpha)}(x)\), where \(j\geq 0\).
The polynomials \(C_j^{(\alpha)}(x)\) are sometimes also called \emph{ultraspherical polynomials} in the literature.

\begin{Satz}
\label{stability-sphere}
The spectral problem \eqref{jacobi-identity-spectral} describing the stability of the identity
map, which we parametrize by \(h_1(x)=-\frac{\pi}{2}+2\arctan(e^x)\), is solved by
\begin{align*}
\xi_j(x)=\frac{1}{\cosh x }C^{(\frac{m+1}{2})}_{j-1}(\tanh x ),\qquad \hat{\lambda}_j=-2m+p+j(j+m-1),
\end{align*}
where \(j\in\,\mathbb{N}_+\).
\end{Satz}

\begin{proof}
In order to solve the eigenvalue problem \eqref{jacobi-identity-spectral} we make the ansatz
$$\xi(x)=\frac{1}{\cosh x }f(x).$$
A direct calculation shows that for this ansatz
\begin{align*}
\frac{d}{dx}\big(\frac{\cosh x}{m}(\xi'(x)-(m-1)\tanh x \xi(x))\big)
=\frac{f''(x)}{m}-\frac{f(x)}{\cosh^2(x)}-\tanh x f'(x)
\end{align*}
such that \eqref{jacobi-identity-spectral} becomes
\begin{align*}
f''(x)-m\tanh x f'(x)+\frac{f(x)}{\cosh^2(x)}\frac{m}{m+p-2}(\hat{\lambda}+m-p)=0.
\end{align*}
We find that \(f(x)=1\) and \(\hat{\lambda}=p-m\) solves the above equation.

In order to completely solve the spectral problem \eqref{jacobi-identity-spectral}
we perform the transformation \(f(x)=u(\tanh x )\) which gives the equation
\begin{align*}
(1-\tanh^2x)u''(\tanh x )-(2+m)\tanh x u'(\tanh x )+\frac{m}{m+p-2}(\hat{\lambda}+m-p)u(\tanh x )=0.
\end{align*}
This equation is of the form \eqref{ode-gegenbauer}, thus the claim follows.

This completes the proof of Theorem \ref{thm:stability-identity}.
\end{proof}

\bibliographystyle{plain}
\bibliography{mybib}

\end{document}